\documentclass[11pt]{article}
\usepackage{amsmath, amsthm, amssymb}
\usepackage{bridges}    
\usepackage{graphicx}
\usepackage[colorlinks=true, urlcolor=blue, citecolor=black, linkcolor=black]{hyperref}

\usepackage{subcaption} 
\usepackage{xcolor}     

\title{Towards Flying through Modular Forms}

\author{David Lowry-Duda\textsuperscript{1} and Adam Sakareassen\textsuperscript{2}
\vspace{10pt}\\
\textsuperscript{1} ICERM \& Brown University, Providence, RI; david@lowryduda.com \\
\textsuperscript{2} {Maths Town YouTube Channel; mathstownchannel@gmail.com}}

\date{}

\begin{document}

\maketitle

\thispagestyle{empty}

\begin{abstract}

Modular forms are highly self-symmetric functions studied in number theory,
with connections to several areas of mathematics. But they are rarely
visualized. We discuss ongoing work to compute and visualize modular forms as
3D surfaces and to use these techniques to make videos flying around
the peaks and canyons of these ``modular terrains.'' Our goal is to make
beautiful visualizations exposing the symmetries of these functions.

\end{abstract}

\section*{Introduction}

We began this project to understand better the shapes of modular forms. A
classical modular form is a complex-valued function, defined naturally either
on the unit disk or on the half-plane $\mathcal{H} = \{ x + i y : x,y \in
\mathbb{R}, y > 0 \}$. Each modular form satisfies a set of functional
equations of the shape
\begin{equation}\label{eq:modular_feq}
  f\left( \frac{az + b}{cz + d} \right) = (cz + d)^k f(z),
\end{equation}
where $(\begin{smallmatrix} a & b \\ c & d \end{smallmatrix})$ is any matrix in a
fixed subgroup $\Gamma \subseteq \mathrm{SL}(2, \mathbb{Z})$ and $k$ is an
integer. These functional equations are restrictive and impose a many
symmetries on the modular form. Each matrix in $\Gamma$ relates $f(z)$ to
another value $f(\frac{az + b}{cz + d})$, but it is hard
to grasp how these relations affect the overall shape of the modular form.

\begin{figure}[h!tbp]
\centering
\begin{minipage}[t]{0.35\textwidth}
	\includegraphics[width=\textwidth]{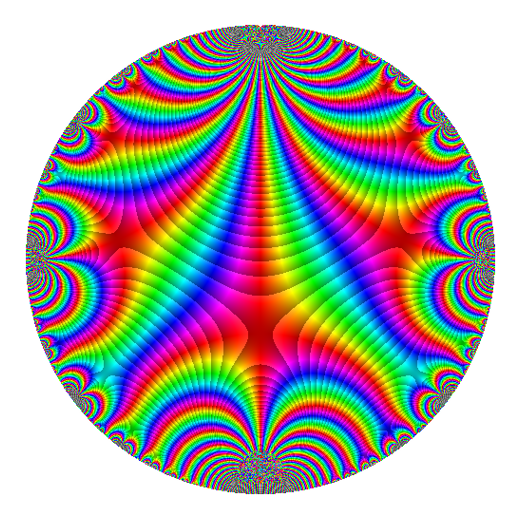}
          \subcaption{%
          }
        	\label{fig:1a}
\end{minipage}
\quad 
\begin{minipage}[t]{0.60\textwidth}
	\includegraphics[width=\textwidth]{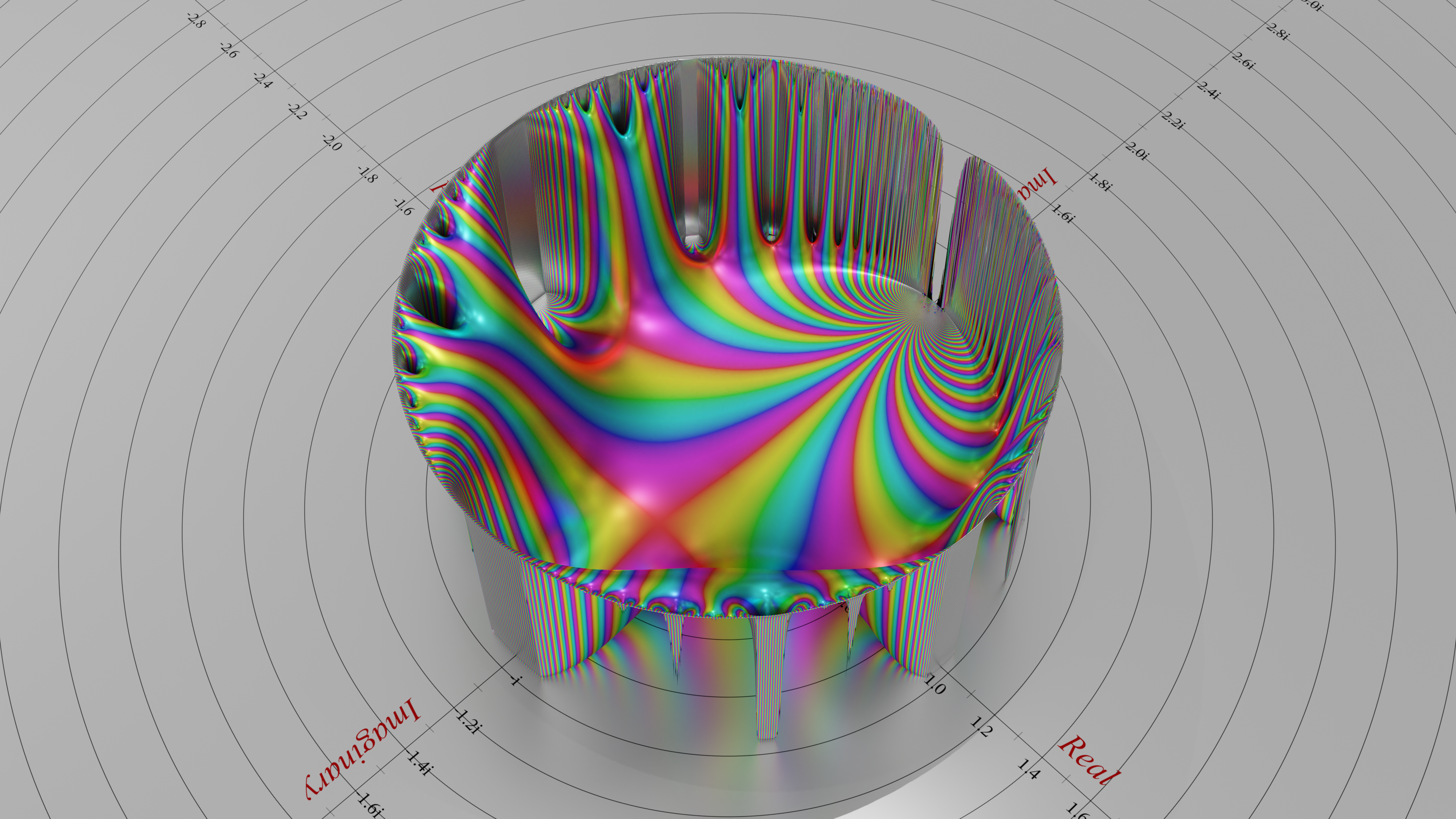}
        	\subcaption{%
          }
        	\label{fig:1b}
\end{minipage}
\caption{A modular form in 2D and 3D.}
\label{fig:1}
\end{figure}

\noindent
In~\cite{LD}, one of the authors began to study various visualizations of
modular forms in detail.
Figure 1(a) 
shows one of the detailed visualizations created. Although a lot of detail is
apparent, we find that this doesn't give a satisfying understanding of the
shape of the modular form. In Figure 1(b), 
we show a particular 3D rendering we created of the same modular form.
In this paper, we describe our ongoing efforts to make accurate, informative,
and beautiful 3D representations similar to Figure 1(b). 

\section*{3D Representation of Modular Forms}

The fundamental problem of visualizing a modular form, either mentally or with
a computer, is that the ``obvious'' visualization is four-dimensional. Thus
when we try to represent a modular form we must carefully choose our
representation. This is true for visualizations of all complex functions $f:
\mathbb{C} \longrightarrow \mathbb{C}$.

Most visualizations for complex functions descend from some form of
\emph{domain coloring}, a term coined initially in Farris's review of
\emph{Visual Complex Analysis}~\cite{FFreview}. To visualize a function $f$
with domain coloring, you first assign a color to each point of the complex
plane. Then you color each point $z$ in the domain of $f$ by the color
assigned to the point $f(z)$. Repeated or predictable choices of coloring often
lead to understandable visualizations, but it is also possible to deliberately
color the plane in a way to emphasize or accentuate behaviors of the visualized
function. We've drawn particular inspiration from~\cite{FF} and~\cite{Wegert},
and hope to explore further color choices in the future.

Our visualizations are 3D domain colorings that we call \emph{modular
terrains}. To each point $z = x + iy$ in the domain, we compute the function
$f(z) = re^{i\theta}$ in polar form. We use the magnitude $r$ to determine the
height of the surface at the point $(x, y)$, and we use the phase $\theta$ to
determine the color of the surface. Finally, we shine a light at the resulting
surface to get a better sense of the topography through highlights and shadows.

Different choices of maps from magnitude to height radically alter the
appearance of the resulting surface. This is particularly noticeable with
modular forms, as the magnitudes will change exponentially, even in small
regions. Figures 1(b), 
2(a), 
and 2(b) 
show three views of the same modular form, each with a different choice of
magnitude-to-height function. In Figure 2(a) 
we used $\tanh(r)$, the hyperbolic tangent function, and then restricted
the output to $[0, 1]$. This separates ``small'' and ``large'' values. In
Figure 2(b) 
we used $\arctan(1/r)$ and then restricted the output to $[0, 1]$.
But we prefer height maps of the form $\arctan(\log(r^\alpha+1))$ for some
$\alpha \in (0, 1)$.
These are the same as the magnitude-to-brightness maps in \S2.2.2 of~\cite{LD}.
In Figure 1(b) 
and later figures, we use $\arctan(\log(\sqrt{r} + 1))$.
The idea behind these maps is that $\log$ tempers the exponential growth of the
magnitudes and $\arctan$ acts like a smooth restriction map, bounding the
heights.
We chose $\sqrt{r}$ through trial and error, as it gives the terrains gentle,
beautiful transitions from lows to highs.

\begin{figure}[h!tbp]
\centering
\begin{minipage}[t]{0.30\textwidth}
	\includegraphics[width=\textwidth]{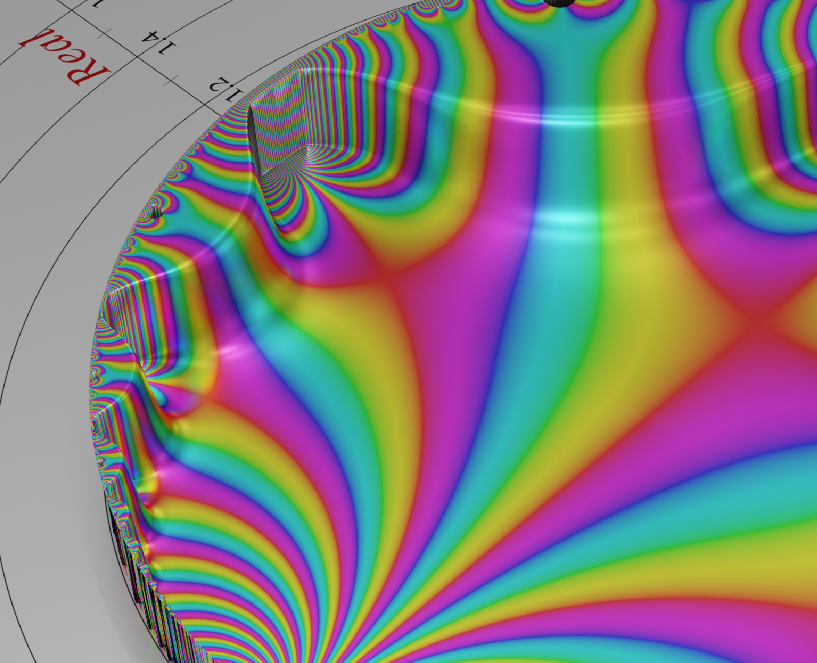}
          \subcaption{%
          }
        	\label{fig:2a}
\end{minipage}
\quad 
\begin{minipage}[t]{0.285\textwidth}
	\includegraphics[width=\textwidth]{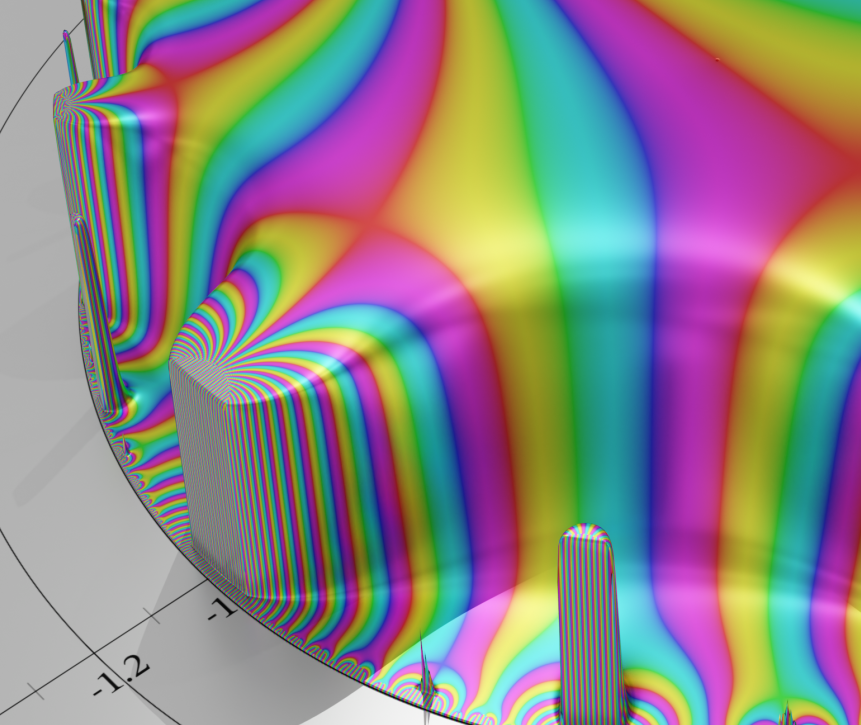}
        	\subcaption{%
          }
        	\label{fig:2b}
\end{minipage}
\caption{The same modular form visualized with two different height functions.}
\label{fig:2}
\end{figure}

We based our color choices on the standard 2D domain coloring,
as in Figure 1(a). 
Color comes from the $\theta$ in $f(z) = r e^{i \theta}$.
We aligned the color wheel so that red corresponds to $\theta = 0$.
But we also project a light onto the modular terrains.
The shadows and highlights affect the colors, but also
serve as visual cues that make the shape more understandable.
Although we sought renderings that clearly describe the modular forms,
we can also use the modular terrains as a canvas.
In Figure 3(d), we ``paint'' on the modular terrain, making art
intrinsically connected to the modular form.

In practice, we compute a triangulation of the surface defined by $f(z)$ and
our choice of height function. We rendered all the images in this paper with
Blender and used shaders for coloring.

\section*{Modular Topography}

Perhaps the most studied modular form is the Delta function\footnote{%
Described further in the LMFDB at
\url{https://www.lmfdb.org/ModularForm/GL2/Q/holomorphic/1/12/a/a/}
}
\begin{equation}\label{eq:delta_exp}
   \Delta(z) = e^{2 \pi i z}\prod_{n \geq 1} (1 - e^{2 \pi i n z})^{24}
   = e^{2 \pi i z} - 24 e^{2 \pi i (2z)} + 252 e^{2 \pi i (3z)} + \cdots
\end{equation}
Figures 1, 2, and 3 all depict $\Delta(z)$. The Delta function connects
work in complex analysis, Riemann surfaces, and algebraic geometry to
the theory of elliptic functions. It is the discriminant of the elliptic
invariants of the Weierstrass $\wp$ function, and is also called the
\emph{modular discriminant}. Ramanujan investigated $\Delta(z)$ and asked
questions leading to the modern theory of modular forms. Our figures depict
$\Delta(z)$.

\begin{figure}[h!tbp]
\centering
\begin{minipage}[t]{0.60\textwidth}
  \begin{minipage}[t]{0.8675\textwidth}
    \includegraphics[width=\textwidth]{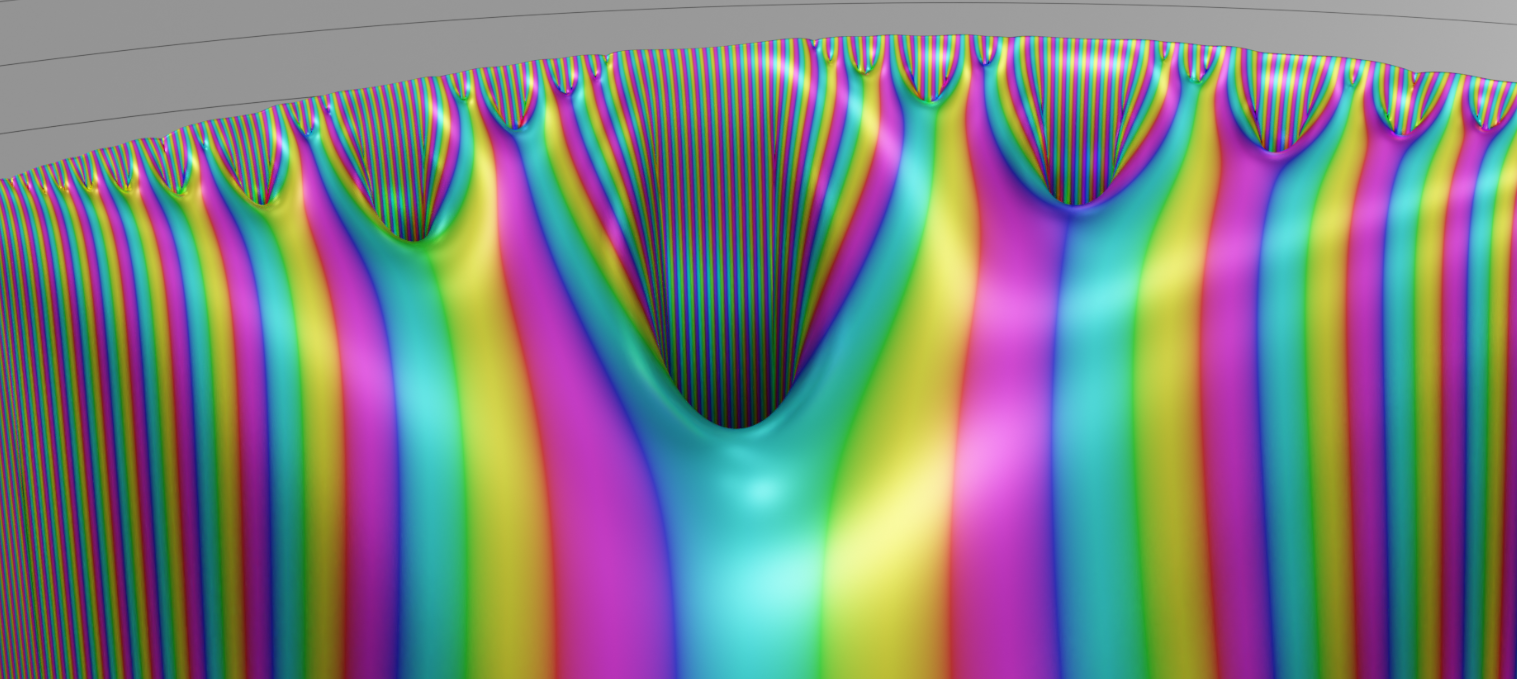}
    \subcaption{}
  \end{minipage}

  \begin{minipage}[t]{0.450\textwidth}
  \includegraphics[width=\textwidth]{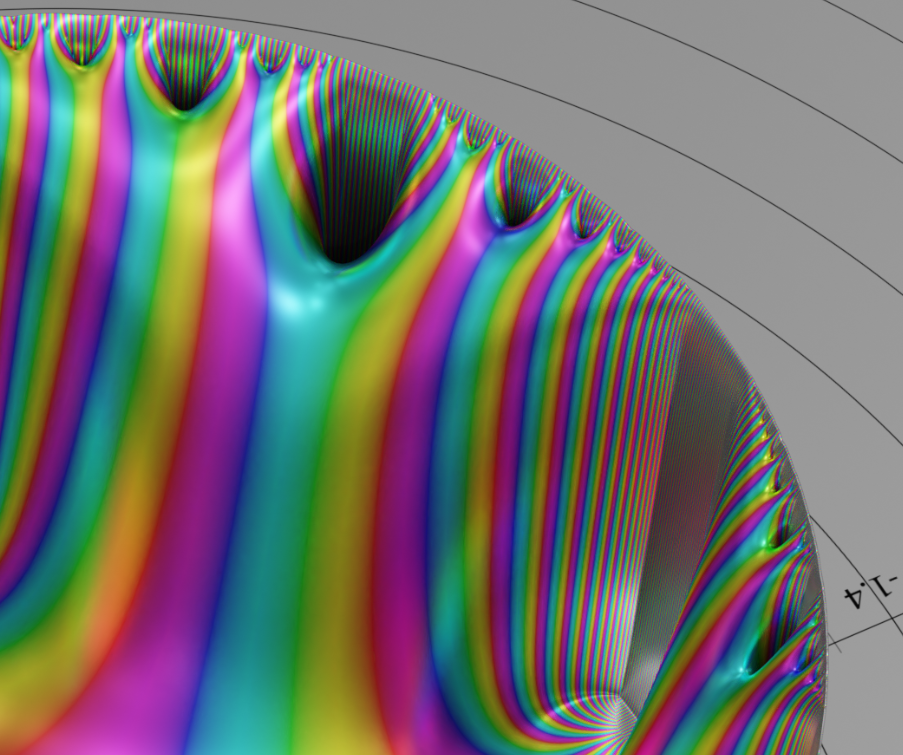}
  \subcaption{}
  \end{minipage}
  ~
  \begin{minipage}[t]{.388\textwidth}
    \includegraphics[width=\textwidth]{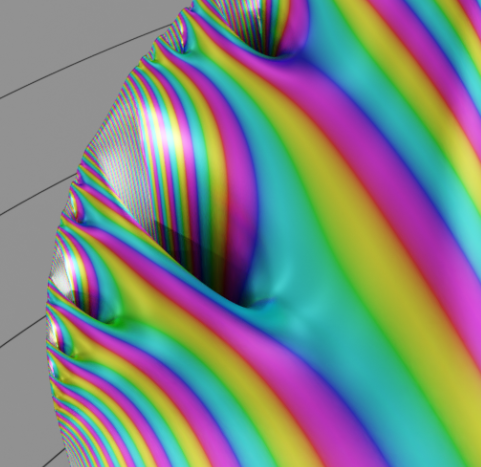}
     \subcaption{}
  \end{minipage}
\end{minipage}
\begin{minipage}{0.35\textwidth}
  \vspace{1.4cm}
  \centering
  \begin{minipage}{\textwidth}
    \includegraphics[width=0.96\textwidth]{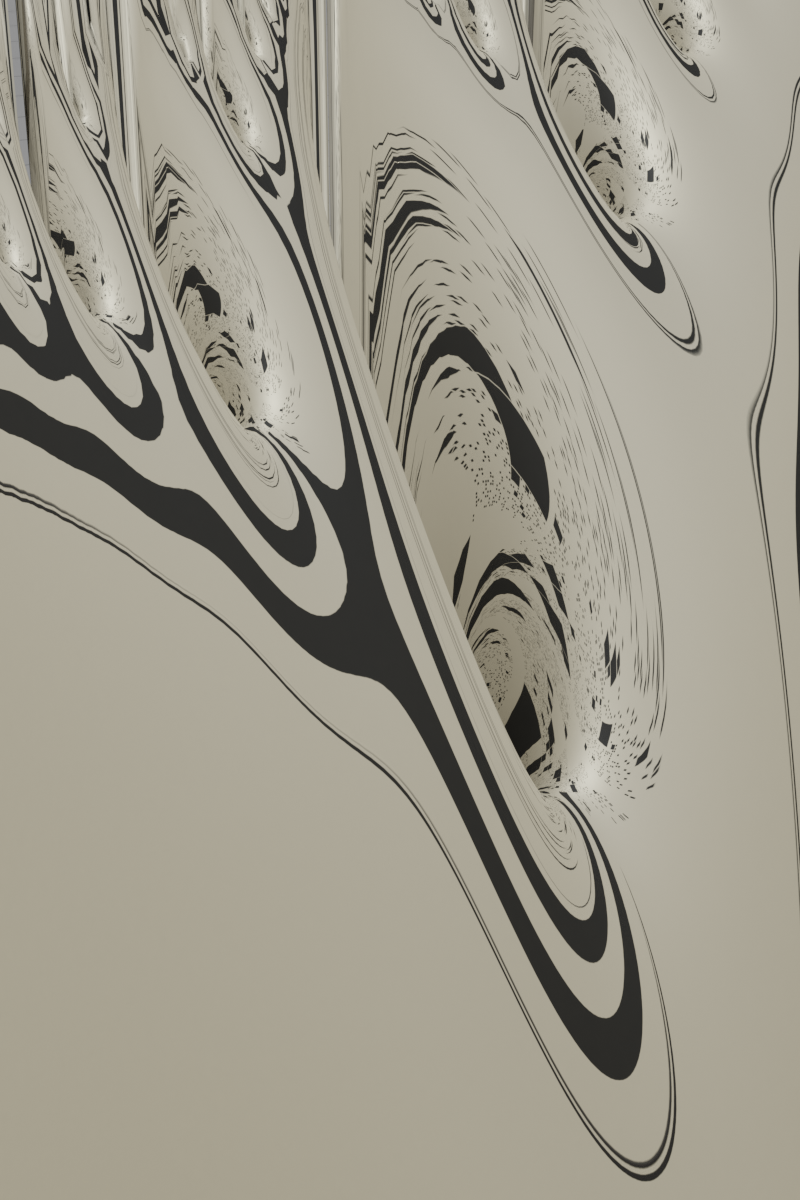}
    \subcaption{}
  \end{minipage}
\end{minipage}
\caption{Four views of the topographic details on the rim of $\Delta(z)$.\\
In (d), we paint the surface with a periodic generative pattern.}
\label{fig:3}
\end{figure}

To produce accurate visualizations, we use the self-symmetries from the
functional equations~\eqref{eq:modular_feq} to compute each value $\Delta(z)$
in a region where the infinite expansion~\eqref{eq:delta_exp} rapidly
converges. In regions where $\Delta(z)$ is extremely small, it is also
necessary to perform the computation with higher machine precision to get a
smooth surface. This enables close examination of the details around a
particular area, even around the ``rim,'' as in Figure~\ref{fig:3}. Notice this
figure gives three points of view of the same region, as well as a beautiful view
with nonstandard coloring. We are working towards producing smooth
video transitions across these points of view, giving the impression of flying
around the modular terrain.

\vspace{-1.1cm}

For other modular forms, we can begin with data in the LMFDB~\cite{LMFDB}
(building on the computations in~\cite{ModularCrew}). We're working on similar
techniques for high precision visualization of arbitrary modular forms, but in
practice the functional equations~\eqref{eq:modular_feq} are more complicated
for different matrix subgroups and precise computation becomes vastly more
tedious near the boundary.

Despite these difficulties, we've begun to compute modular terrains for
different modular forms. In Figure~4(a) 
we see the modular form \texttt{5.4.a.a}\footnote{Corresponding LMFDB page:
\url{https://www.lmfdb.org/ModularForm/GL2/Q/holomorphic/5/4/a/a/}}, and in
Figure~4(b) 
we see the modular form \texttt{56.1.h.a}\footnote{Corresponding LMFDB page:
\url{https://beta.lmfdb.org/ModularForm/GL2/Q/holomorphic/56/1/h/a/}}. We note
that the number and variety of ``canyons'' visually differ, as does the rate of
``elevation change'' --- and further, we find this much clearer in 3D than in
2D. We look forward to visualizing more modular forms.

\begin{figure}[h!tbp]
\centering
\begin{minipage}[t]{0.365\textwidth}
\includegraphics[width=\textwidth]{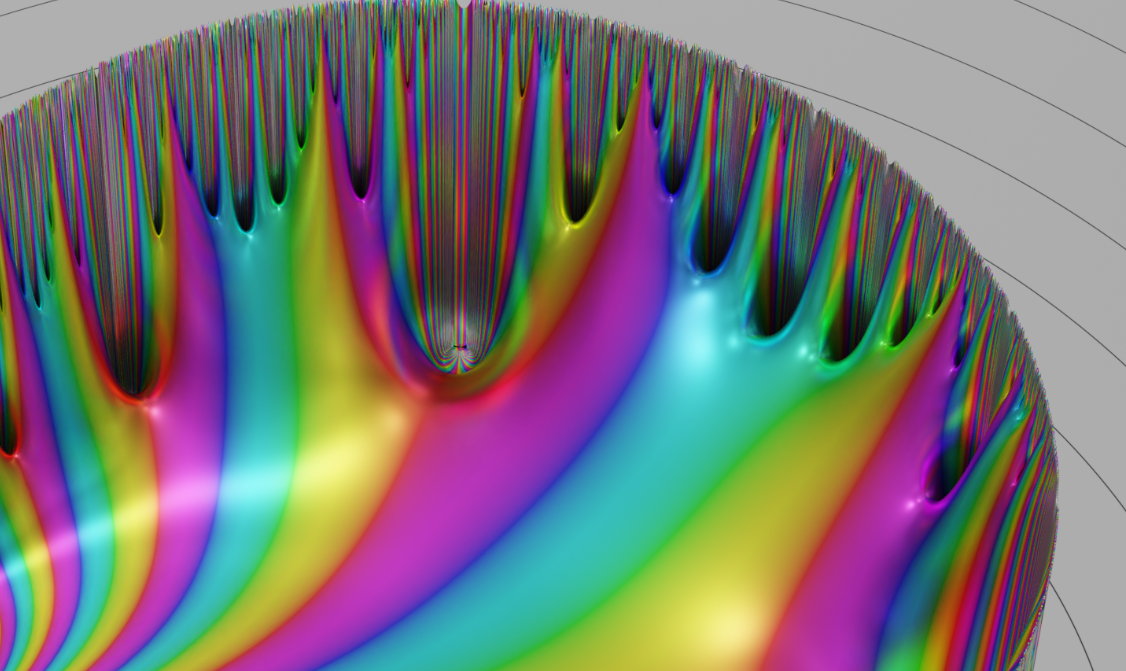}
\subcaption{}\label{fig:4a}
\end{minipage}
~
\begin{minipage}[t]{.365\textwidth}
	\includegraphics[width=\textwidth]{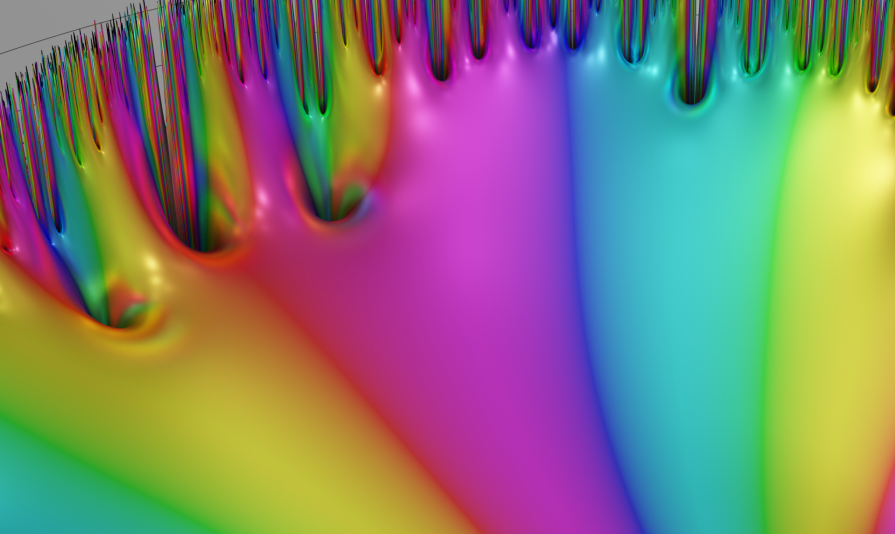}
   \subcaption{}\label{fig:4b}
\end{minipage}
\caption{Detail views of the modular forms \texttt{5.4.a.a} (a) and
\texttt{56.1.h.a} (b).}
\label{fig:4}
\end{figure}

\vspace{-0.75cm}

\section*{Summary and Future Work}

We've described a set of 3D visualizations and presented several static
renderings.
In the YouTube video \url{https://www.youtube.com/watch?v=s6sdEbGNdic}, we have
our first \emph{dynamic}, flying view of $\Delta(z)$.
We're working on efficient methods of computing visualizations of
more modular forms, which we plan on using to make more
visualizations. These will appear on the second author's YouTube
channel ``The Mathemagicians' Guild,''
\url{https://www.youtube.com/channel/UCHsYWDqJkNlboBpyrRqLzkA}.
We also plan on making more art with these terrains as canvas, as in
Figure 3(d).

\vspace{-0.25cm}
\section*{Acknowledgements}

DLD was supported by the Simons Collaboration in Arithmetic Geometry, Number
Theory, and Computation via the Simons Foundation grant 546235. This project
owes much to many conversations at the fall 2019 ICERM program on Illustrating
Mathematics.

\vspace{-0.25cm}
{\setlength{\baselineskip}{13pt}
\raggedright

} 

\end{document}